\documentclass[11pt, a4paper]{amsart}
\usepackage[utf8]{inputenc}
\usepackage{amssymb,amsthm,amsmath, amsfonts}

\newcommand{\mb}[1]{\mathbb{{#1}}}
\newcommand{\e}{\varepsilon}
\newcommand{\call}[4]{\int_{#1}^{#2} {#3} \; \textrm{d} {#4}}
\newcommand{\dd}{\textrm{d}}

\newtheorem*{que*}{Question}

\theoremstyle{remark}
\newtheorem*{remark*}{Remark}

\title[A Note on a Brunn-Minkowski Inequality for the Gaussian Measure]{A Note on a Brunn-Minkowski Inequality for the Gaussian Measure}
\author{Piotr Nayar}
\address{Piotr Nayar \\ Institute of Mathematics, University of Warsaw, Banacha 2, 02-097 Warszawa, Poland.}
\email{nayar@mimuw.edu.pl}
\thanks{Research of the first author partially supported by NCN Grant no. 2011/01/N/ST1/01839.}

\author{Tomasz Tkocz}
\address{Tomasz Tkocz \\ Institute of Mathematics, University of Warsaw, Banacha 2, 02-097 Warszawa, Poland.}
\email{t.tkocz@mimuw.edu.pl}
\thanks{Research of the second author partially supported by NCN Grant no. 2011/01/N/ST1/05960.}

\keywords{Convex body, Gauss measure, Brunn-Minkowski inequality, B-conjecture}
\subjclass[2010]{Primary 52A40; Secondary 60G15}

\begin{document}

\begin{abstract}
  We give the counter-examples related to a Gaussian Brunn-Minkowski inequality and the (B) conjecture.
\end{abstract}

\maketitle

\section{Introduction and notation}\label{sec:intro}

Let $\gamma_n$ be the standard Gaussian distribution on $\mb{R}^n$, i.e. the measure with the density
\[
  g_n(x)=\frac{1}{(2\pi)^{n/2}} e^{-|x|^2/2},
\]
where $|\cdot |$ stands for the standard Euclidean norm. A powerful tool in convex geometry is the Brunn-Minkowski inequality for Lebesgue measure (see \cite{Sch} for more information). Concerning the Gaussian measure, the following question has recently been posed.
\begin{que*}[R. Gardner and A. Zvavitch, \cite{GZ}]
Let $0 < \lambda < 1$ and let $A$ and $B$ be closed convex sets in $\mb{R}^n$ such that $o \in A \cap B$. Is it true that   
\begin{equation}\label{eq.GBM}\tag{GBM}
  \gamma_n(\lambda A + (1-\lambda)B)^{1/n} \geq \lambda \gamma_n(A)^{1/n} + (1-\lambda) \gamma_n(B)^{1/n}?
\end{equation}  
\end{que*}
A counter-example is given in this note. However, we believe that this question has an affirmative answer in the case of $o$-symmetric convex sets, i.e. the sets satisfying $K = -K$.

In \cite{CFM} it is proved that for an $o$-symmetric convex set $K$ in $\mb{R}^n$ the function
\begin{equation}\label{eq.fun}
  \mb{R} \ni t \ \mapsto \gamma_n(e^tK),
\end{equation}
is log-concave. This was conjectured by W. Banaszczyk and popularized by R. Lata\l a \cite{Lat}. It turns out that the (B) conjecture cannot be extended to the class of sets which are not necessarily $o$-symmetric yet contain the origin, as one of the sets provided in our counter-example shows.

As for the notation, we frequently use the function 
\[
  T(x) = \frac{1}{\sqrt{2\pi}}\int_x^\infty e^{-t^2/2} \dd t.
\]

\section{Counter-examples}\label{sec.counter-examples}

Now we construct the convex sets $A, B \subset \mb{R}^2$ containing the origin such that inequality \eqref{eq.GBM} does not hold. Later on we show that for the set $B$ the (B) conjecture is not true.

Fix $\alpha \in (0,\pi/2)$ and $\e>0$. Take
\begin{align*}
	A &= \{ (x,y) \in \mb{R}^2 \; | \; y \geq |x| \tan \alpha \}, \\
	B = B_\e &= \{ (x,y) \in \mb{R}^2 \; | \; y \geq |x| \tan \alpha - \e  \} = A -(0,\e).
\end{align*}
Clearly, $A,B$ are convex and $0\in A \cap B$. Moreover, from convexity of $A$ we have $\lambda A + (1-\lambda)A=A$ and therefore
\[
	\lambda A + (1-\lambda)B = \lambda A + (1-\lambda)(A -(0,\e)) = A - (1-\lambda)(0,\e).
\]
Observe that
\begin{align*}
\gamma_2(A) &= \frac{1}{2} - \frac{\alpha}{\pi}, \\
\gamma_2(B) &= 2 \call{0}{+\infty}{T(x \tan \alpha - \e) \frac{1}{\sqrt{2\pi}}e^{-x^2/2}}{x}, \\
\gamma_2(\lambda A + (1-\lambda)B) &= 2 \call{0}{+\infty}{T(x \tan \alpha - \e(1-\lambda)) \frac{1}{\sqrt{2\pi}}e^{-x^2/2}}{x}
\end{align*}
and that these expressions are analytic functions of $\e$. We will expand these functions in $\e$ up to the order $2$. Let 
\[
	a_k = \call{0}{+\infty}{T^{(k)}(x \tan \alpha)\frac{1}{\sqrt{2\pi}}e^{-x^2/2}}{x},
\]
for $k = 0,1,2$, where $T^{(k)}$ is the $k$-th derivative of $T$ (we adopt the standard notation $T^{(0)}=T$).
We get
\begin{align*}
\gamma_2(A) &= 2a_0, \\
\gamma_2(B) &= 2a_0 - 2\e a_1 + \e^2 a_2 + o(\e^2), \\
\gamma_2(\lambda A + (1-\lambda)B) &= 2a_0 - 2\e(1-\lambda) a_1 + \e^2 (1-\lambda)^2 a_2 + o(\e^2).
\end{align*}
Thus
\[
\sqrt{\gamma_2(B)} = \sqrt{2a_0} - \frac{a_1}{\sqrt{2a_0}} \e + \left( \frac{a_2}{2\sqrt{2a_0}} - \frac{a_1^2}{2 (2a_0)^{3/2}}   \right)\e^2 + o(\e^2). 
\]
Taking $\e(1-\lambda)$ instead of $\e$ we obtain 
\begin{align*}
\sqrt{ \gamma_2(\lambda A + (1-\lambda)B)} & = \sqrt{2a_0} - \frac{a_1}{\sqrt{2a_0}}(1-\lambda) \e \\ & \quad +  \left( \frac{a_2}{2\sqrt{2a_0}} 
 \quad - \frac{a_1^2}{2 (2a_0)^{3/2}}   \right)(1-\lambda)^2\e^2 + o(\e^2).  
\end{align*}
Since
\begin{gather*} \sqrt{ \gamma_2(\lambda A + (1-\lambda)B)} - \lambda \sqrt{\gamma_2(A)} - (1-\lambda)\sqrt{\gamma_2(B)} \\
= -\lambda (1-\lambda) \frac{1}{2 (2a_0)^{3/2}}(2a_0 a_2 - a_1^2)\e^2 + o(\e^2),
\end{gather*}
we will have a counter-example if we find $\alpha \in (0,\pi/2)$ such that 
\[
	2a_0 a_2 - a_1^2 > 0.
\] 
Recall that $a_0 = \frac{1}{2}\gamma_2(A) = \frac{1}{2}\left( \frac{1}{2} - \frac{\alpha}{\pi} \right)$. The integrals that define the $a_k$'s can be calculated. Namely,
\begin{align*}
   a_1 &= \call{0}{\infty}{ T'(x\tan \alpha) \frac{e^{-x^2/2}}{\sqrt{2\pi}} }{x} = -\frac{1}{\sqrt{2\pi}}\frac{1}{2}\int_{\mb{R}} e^{-(1+\tan^2 \alpha)x^2/2} \frac{\dd x}{\sqrt{2\pi}} \\
   &= -\frac{1}{\sqrt{2\pi}}\frac{1}{2\sqrt{1+\tan^2 \alpha}}, \\   
   a_2 &= \call{0}{\infty}{ T''(x\tan \alpha) \frac{e^{-x^2/2}}{\sqrt{2\pi}} }{x} = \frac{1}{\sqrt{2\pi}}\int_{0}^{\infty} (x\tan \alpha) e^{-(1+\tan^2\alpha)x^2/2} \frac{\dd x}{\sqrt{2\pi}} \\
   &= \frac{1}{2\pi}\frac{\tan \alpha}{1+\tan^2 \alpha}.
\end{align*}
Therefore,
\begin{align*} 2a_0 a_2 -a_1^2 & = 2 \left( \frac12 \left(\frac12 - \frac{\alpha}{\pi} \right) \cdot \frac{1}{2\pi} \frac{\tan \alpha}{1+ \tan^2 \alpha} \right) - \frac{1}{2\pi} \cdot \frac{1}{4(1+\tan^2 \alpha)} \\
& = \frac{1}{8\pi} \frac{1}{1+\tan^2 \alpha} \left( \tan \alpha \left( 2 - \frac{4\alpha}{\pi} \right) - 1  \right),
\end{align*}
which is positive for $\alpha$ close to $\pi/2$.

Now we turn our attention to the (B) conjecture. We are to check that for the set $B=B_\e$ the function $\mb{R} \ni t \ \mapsto \gamma_n(e^tB)$ is not log-concave, provided that $\e$ is sufficiently small. Since
\[
  e^t B = \{(x, y) \in \mb{R}^2 \ | \ y \geq \tan \alpha |x| - \e e^t \}
\]
we get
\begin{align*}
  \ln \gamma_2(e^tB) &= \ln \Bigg( 2\int_0^\infty T(x\tan \alpha - e^t\e )\frac{e^{-x^2/2}}{\sqrt{2\pi}}\dd x \Bigg) \\
  &= \ln\left( 2\int_0^\infty T(x\tan \alpha)\frac{e^{-x^2/2}}{\sqrt{2\pi}}\dd x \right) - \e e^t\frac{\int_0^\infty T'(x\tan \alpha)e^{-x^2/2} \dd x}{\int_0^\infty T(x\tan \alpha)e^{-x^2/2} \dd x} + o(\e).
\end{align*}
This produces the desired counter-example for sufficiently small $\e$ as the function $t \mapsto \beta e^t$, where
\[
  \beta = -\frac{\int_0^\infty T'(x\tan \alpha)e^{-x^2/2} \dd x}{\int_0^\infty T(x\tan \alpha)e^{-x^2/2} \dd x} > 0,
\]
is convex.\qed

\begin{remark*}
The set $B_\e$ which serves as a counter-example to the (B) conjecture in the nonsymmetric case works when the parameter $\alpha = 0$ as well (and $\e$ is sufficiently small). Since $B_\e$ is simply a halfspace in this case, it shows that symmetry of $K$ is required for log-concavity of \eqref{eq.fun} even in the one-dimensional case.
\end{remark*}

\section*{Acknowledgements}

The authors would like to thank Professors R. Gardner and A. Zvavitch for pointing out that the constructed set may also serve as a counter-example to the (B) conjecture in the non-symmetric case. An anonymous referee deserves thanks for the remark.

\end{document}